\documentclass[11pt,reqno]{amsart}

\usepackage[T1]{fontenc}
\usepackage[utf8]{inputenc}
\usepackage{mathpazo}
\usepackage[scaled=.95]{helvet}
\usepackage{microtype}
\usepackage[margin=1.10in]{geometry}
\usepackage{amsmath,amssymb,mathtools,mathrsfs}
\usepackage{xcolor}
\usepackage{enumitem}
\usepackage{hyperref}
\usepackage{aliascnt}
\usepackage[nameinlink,noabbrev]{cleveref}

\hypersetup{
  colorlinks=true,
  linkcolor=blue!55!black,
  citecolor=blue!55!black,
  urlcolor=blue!55!black,
  pdfauthor={Juan J. Manfredi, Aelson Sobral, Eduardo V. Teixeira, and Jos\'e Miguel Urbano},
  pdftitle={Higher-order estimates for p-harmonic functions in the plane}
}
\numberwithin{equation}{section}
\theoremstyle{plain}
\newtheorem{theorem}{Theorem}[section]

\newaliascnt{proposition}{theorem}
\newtheorem{proposition}[proposition]{Proposition}
\aliascntresetthe{proposition}

\newaliascnt{lemma}{theorem}
\newtheorem{lemma}[lemma]{Lemma}
\aliascntresetthe{lemma}

\newaliascnt{corollary}{theorem}
\newtheorem{corollary}[corollary]{Corollary}
\aliascntresetthe{corollary}

\theoremstyle{remark}
\newaliascnt{remark}{theorem}
\newtheorem{remark}[remark]{Remark}
\aliascntresetthe{remark}

\newcommand{\R}{\mathbb R}
\newcommand{\N}{\mathbb N}
\newcommand{\norm}[1]{\left\lVert #1\right\rVert}
\newcommand{\abs}[1]{\left\lvert #1\right\rvert}
\newcommand{\diver}{\operatorname{div}}
\newcommand{\Kcal}{\mathcal K}

\setlist[enumerate]{leftmargin=2.2em,itemsep=4pt,topsep=4pt}
\setlist[itemize]{leftmargin=2em,itemsep=3pt,topsep=4pt}

\title[Higher-order estimates in the plane]{Higher-order estimates for $p$-harmonic functions in the plane}

\author[J.J. Manfredi]{Juan J. Manfredi}
\address{Department of Mathematics, University of Pittsburgh, 301 Thackeray Hall, Pittsburgh, PA 15260, USA}
\email{manfredi@pitt.edu}

\author[A. Sobral]{Aelson Sobral}
\address{Applied Mathematics and Computational Sciences (AMCS), Computer, Electrical and Mathematical Sciences and Engineering Division (CEMSE), King Abdullah University of Science and Technology (KAUST), Thuwal, 23955-6900, Kingdom of Saudi Arabia}
\email{aelson.sobral@kaust.edu.sa}

\author[E. V. Teixeira]{Eduardo V. Teixeira}
\address{Department of Mathematics, Oklahoma State University, 401 MSCS, Stillwater, OK 74078-1058, USA}
\email{eduardo.teixeira@okstate.edu}

\author[J.M.~ Urbano]{Jos\'{e} Miguel Urbano}
\address{Applied Mathematics and Computational Sciences (AMCS), Computer, Electrical and Mathematical Sciences and Engineering Division (CEMSE), King Abdullah University of Science and Technology (KAUST), Thuwal, 23955-6900, Kingdom of Saudi Arabia; and CMUC, Department of Mathematics, University of Coimbra, 3000-143 Coimbra, Portugal}
\email{miguel.urbano@kaust.edu.sa} 

\date{}

\begin{document}

\begin{abstract}
Closing a longstanding gap in the theory, we establish uniform interior $C^{k,\beta}$ estimates for planar $p$-harmonic functions with $1<p<2$. These estimates hold for every integer $k\ge2$ and every $0<\beta<1$ such that $k+\beta<s(p)$, where $s(p)$ is the optimal regularity order identified by Iwaniec and Manfredi. As consequences, we obtain decay estimates at critical points, a Liouville theorem with optimal growth exponent, and uniform higher-order estimates for the flux of planar $q$-harmonic functions with $q>2$.
\end{abstract} 

\maketitle

\section{Introduction and main theorem}

We establish uniform higher-order interior estimates for weak solutions of
\begin{equation}\label{eq:plap}
	-\Delta_pu\coloneqq-\diver\bigl(\abs{Du}^{p-2}Du\bigr)=0 \quad\text{in }B_1,\qquad 1<p<2.
\end{equation}

Solutions of \eqref{eq:plap} satisfy interior $C^{1,\kappa}$ estimates for some $\kappa\in(0,1)$ depending only on $p$; see Lewis \cite{Lewis1983} and Tolksdorf \cite{Tolksdorf1984}. In the plane, quasiregular methods give more precise regularity information, as developed by Bojarski and Iwaniec \cite{BojarskiIwaniec1987} and Manfredi \cite{Manfredi1988}. Iwaniec and Manfredi \cite{IwaniecManfredi1989} subsequently determined the optimal qualitative regularity (see also Aronsson \cite{Aronsson1989}), with total order
\begin{equation}\label{eq:IM}
	s(p)\coloneqq\frac16\left(7+\frac1{p-1}+\sqrt{1+\frac{14}{p-1}+\frac1{(p-1)^2}}\right).
\end{equation}
Their theorem gives local $C^{m,\alpha}$ regularity, where $m\in\N$, $0<\alpha\le1$, and $m+\alpha=s(p)$. For $1<p<2$, this order satisfies $s(p)>2$ and tends to infinity as $p\downarrow1$.

The qualitative theorem does not, by itself, provide uniform bounds for the corresponding norms. Danielli and Petrosyan \cite[Remark~2.6]{DanielliPetrosyan2006} already noted the difficulty of controlling the constants in this higher regularity theory by the size of the solution. Ara\'ujo, Teixeira, and Urbano discuss the same distinction for the sharp planar gradient exponent \cite[Section~2]{AraujoTeixeiraUrbano2017}, where the quantitative input comes from the work of Baernstein and Kovalev \cite{BaernsteinKovalev2005} on quasiregular gradient mappings; see also the remarks in \cite[Section~2]{BalciDieningWeimar2020}. Thus, obtaining estimates controlled only by a solution norm requires an additional argument.

Related quantitative results address gradient regularity and behavior at critical points. Balci, Diening, and Weimar \cite[Theorem~2.3]{BalciDieningWeimar2020} proved gradient oscillation estimates with every H\"older exponent below one for planar $p$-harmonic functions with $1<p\le2$, as part of their work on higher-order estimates for the nonlinear flux. Teixeira \cite{Teixeira2022} obtained improved regularity estimates at critical points of $p$-degenerate equations in heterogeneous media when the source vanishes there, with the gains governed by its rate of vanishing; see also \cite{Teixeira2014}. Here we establish uniform estimates for the higher classical derivatives of planar $p$-harmonic functions throughout the range of orders strictly below $s(p)$.

\subsection{Main result}

Throughout the paper, $B_r=B_r(0)\subset\R^2$ and $1<p<2$. A function $u\in W^{1,p}_{\mathrm{loc}}(B_1)$ is $p$-harmonic in $B_1$ if
\begin{equation*}
	\int_{B_1}\abs{Du}^{p-2}Du\cdot D\varphi\,dx=0 \qquad\text{for every }\varphi\in C_c^\infty(B_1).
\end{equation*}
Our main result is the following.

\begin{theorem}\label{thm:main}
Let $1<p<2$. For every integer $k\ge2$ and every $0<\beta<1$ satisfying $k+\beta<s(p)$, there is $C=C(p,k,\beta)<\infty$ such that every $p$-harmonic function $u$ in $B_1$ with $u\in L^p(B_1)$ satisfies
\begin{equation}\label{eq:main}
	\sum_{j=0}^k\norm{D^ju}_{L^\infty(B_{1/2})} +[D^ku]_{C^{0,\beta}(B_{1/2})} \le C\norm{u}_{L^p(B_1)}.
\end{equation}
\end{theorem}

To our knowledge, \Cref{thm:main} gives the first uniform interior $C^{k,\beta}$ estimates for planar $p$-harmonic functions in the singular range $1<p<2$ that cover every $k\ge2$ and $0<\beta<1$ with $k+\beta<s(p)$. These bounds depend only on $p,k,\beta$ and the $L^p$ norm of the solution.

If $s(p)=m+\alpha$ with $m\in\N$ and $0<\alpha<1$, the theorem gives uniform $C^{m,\beta}$ estimates for every $0<\beta<\alpha$. If $s(p)=m$ is an integer, it gives uniform $C^{m-1,\beta}$ estimates for every $0<\beta<1$. Since $s$ is continuous and strictly decreasing on $(1,2)$, with $s(p)\to\infty$ as $p\downarrow1$ and $s(p)\to2$ as $p\uparrow2$, the second case occurs for exactly one exponent $p_m$ for each integer $m\ge3$; solving $s(p)=m$ in \eqref{eq:IM} gives
\begin{equation*}
	p_m=1+\frac{m}{(3m-4)(m-1)},\qquad\text{that is,}\qquad p_3=\tfrac{13}{10},\quad p_4=\tfrac76,\quad p_5=\tfrac{49}{44},\ \ldots
\end{equation*}

\begin{remark}\label{rem:scope}
\begin{enumerate}[label=(\roman*)]
\item The range $k+\beta<s(p)$ cannot be extended beyond $s(p)$, because the regularity order of Iwaniec and Manfredi is optimal \cite{IwaniecManfredi1989}. We do not know whether a uniform estimate holds at the endpoint order $s(p)$ itself. The compactness argument below uses the strict inequality in an essential way (see the proof of \Cref{lem:scales}) and gives no explicit information on the constant $C(p,k,\beta)$.
\item The norm $\norm u_{L^p(B_1)}$ in \eqref{eq:main} can be replaced by $\norm u_{L^q(B_1)}$ for any $q>p-1$, with $C$ depending also on $q$; see \cite[Lemma~3.6]{Lindqvist}.
\end{enumerate}
\end{remark}

\subsection{Consequences}

\Cref{thm:main} and its proof have three consequences, which we collect in \Cref{sec:consequences}. First, at a critical point $x_0$ the difference $u(x)-u(x_0)$ decays at every rate below $s(p)$, with a uniform constant:
\begin{equation*}
	\abs{u(x)-u(x_0)}\le C\norm u_{L^p(B_1)}\abs{x-x_0}^\sigma\qquad\text{for every }\sigma<s(p);
\end{equation*}
see \Cref{cor:critical}. Second, a $p$-harmonic function in the whole plane that grows at most as $\abs x^\gamma$ with $\gamma<s(p)$ is affine, and the exponent $s(p)$ is optimal for this conclusion; see \Cref{cor:liouville} and \Cref{rem:sharp}. Third, the stream function of a planar $q$-harmonic function is $q'$-harmonic, with $q'=q/(q-1)$. For $q>2$, \Cref{thm:main} therefore gives uniform $C^{j,\beta}$ estimates for the flux $\abs{Dv}^{q-2}Dv$ of a $q$-harmonic function $v$, for every integer $j\ge1$ and $0<\beta<1$ with
\begin{equation*}
	j+\beta<\tfrac16\left(q+\sqrt{q^2+12q-12}\right);
\end{equation*}
see \Cref{cor:flux}. This complements the first-order flux estimates of \cite[Theorem~2.2]{BalciDieningWeimar2020}.

\subsection{Contribution and proof strategy}

We follow the qualitative-to-quantitative approach of Sobral, Teixeira, and Urbano \cite{SobralTeixeiraUrbano}. Their $p$-Laplace argument for $p>2$ yields gradient estimates through affine approximation and compactness for shifted equations. The additional ingredient here is a consequence of Lewis's planar polynomial rigidity theorem \cite{Lewis1980}: for $p\ne2$, a homogeneous $p$-harmonic polynomial of degree at least two in the plane vanishes identically. By a blow-up argument, it follows that at a critical point, every derivative through order $k$ vanishes whenever the solution is $C^k$. Together with Iwaniec--Manfredi regularity, this gives decay of order strictly greater than $k+\beta$ after subtracting the value at that point.

Compactness turns this decay into a uniform improvement at one of finitely many scales for normalized solutions with small central gradient. We rescale and repeat as long as the gradient remains small. If the iteration continues indefinitely, the Taylor polynomial is constant. If it stops, classical estimates on a ball where the gradient stays away from zero give a Taylor polynomial of degree at most $k$. Scaling then controls its coefficients and the approximation error at all larger radii. Comparing the resulting Taylor polynomials at nearby centers yields the H\"older estimate for $D^ku$.

We establish the preliminary estimates in \Cref{sec:critical}, prove the main theorem in \Cref{sec:estimates}, and derive the consequences in \Cref{sec:consequences}.

\section{Critical points and compactness}\label{sec:critical}

Fix $1<p<2$ and an integer $k\ge2$. We establish the vanishing of higher-order derivatives at critical points, recall classical estimates at noncritical points, and use compactness to obtain a uniform decay estimate at one of finitely many scales.

\begin{lemma}\label{lem:rigidity}
Let $u$ be $p$-harmonic and $C^k$ in a neighborhood of a point $x_0$ with $Du(x_0)=0$. Then
\begin{equation}\label{eq:vanishing}
	D^ju(x_0)=0\quad \text{for} \quad  j \in \{2,\ldots,k\}.
\end{equation}
\end{lemma}

\begin{proof}
We argue by induction on $m=2,\ldots,k$. Assume that $D^ju(x_0)=0$ for $1\le j<m$. Define
\begin{equation*}
	P_m(x)=\frac{1}{m!}D^mu(x_0)[x^m],\qquad u_r(x)=\frac{u(x_0+rx)-u(x_0)}{r^m},
\end{equation*}
where $D^mu(x_0)[x^m]$ denotes evaluation of the derivative tensor on $m$ copies of $x$. Taylor's theorem applied to $u$ and $Du$ gives $u_r\to P_m$ in $C^1_{\mathrm{loc}}(\R^2)$ as $r\downarrow0$. Since each $u_r$ is $p$-harmonic, passing to the limit in the weak equation shows that $P_m$ is also $p$-harmonic. Lewis's planar polynomial rigidity theorem \cite[Lemma~2]{Lewis1980} therefore yields $P_m\equiv0$; see also the introduction of \cite{LewisVogel2016}. Hence $D^mu(x_0)=0$, completing the induction.
\end{proof}

\subsection{Estimates at noncritical points}

The classical first-order estimates \cite{Lewis1983,Tolksdorf1984} give $\kappa\in(0,1)$ and $C_0$, depending only on $p$, such that
\begin{equation}\label{eq:first}
	\norm{Dv}_{L^\infty(B_{3/4})} +[Dv]_{C^{0,\kappa}(B_{3/4})}\le C_0
\end{equation}
whenever $v$ is $p$-harmonic in $B_1$ and $\norm v_{L^\infty(B_1)}\le1$. We recall the following Schauder estimates.

\begin{lemma}\label{lem:noncritical}
For every $\delta>0$, there are $\tau\in(0,1/4)$ and $C_\delta<\infty$, depending only on $p,k,\delta$, such that every function $v$ satisfying
\begin{equation}\label{eq:norm}
	\Delta_pv=0\text{ in }B_1,\qquad \norm v_{L^\infty(B_1)}\le1,\qquad \abs{Dv(0)}>\delta
\end{equation}
obeys
\begin{equation}\label{eq:noncritical}
	\norm v_{C^{k+1}(B_\tau)}\le C_\delta.
\end{equation}
\end{lemma}

\begin{proof}
Choose $\tau\in(0,1/4)$ so that $C_0(2\tau)^\kappa\le\delta/2$. By \eqref{eq:first}, $\abs{Dv(x)-Dv(0)}\le\delta/2$ and hence
\begin{equation*}
	\delta/2\le\abs{Dv(x)}\le C_0\qquad\text{in }B_{2\tau}.
\end{equation*}
Set $A(\xi)=\abs{\xi}^{p-2}\xi$. On the annulus $\delta/2\le\abs{\xi}\le C_0$, the map $A$ is smooth with bounded derivatives, and the matrix
\begin{equation*}
	D_\xi A(\xi)=\abs{\xi}^{p-2}\left(I+(p-2)\frac{\xi\otimes\xi}{\abs{\xi}^2}\right)
\end{equation*}
is uniformly elliptic, with ellipticity bounds depending only on $p,\delta$.

Since $Dv$ does not vanish in $B_{2\tau}$, the function $v$ is smooth there; see \cite[Section~4]{Lindqvist}. Differentiating the equation gives
\begin{equation*}
	\diver\bigl(D_\xi A(Dv)D(v_{x_i})\bigr)=0\quad\text{in }B_{2\tau},\qquad i=1,2.
\end{equation*}
By \eqref{eq:first}, the coefficients have uniformly controlled $C^{0,\kappa}$ norms. Interior divergence-form Schauder estimates \cite[Theorem~8.32]{GilbargTrudinger} therefore give $C^{1,\kappa}$ bounds for $v_{x_i}$, and hence $C^{2,\kappa}$ bounds for $v$, on smaller balls. The coefficients then have controlled $C^{1,\kappa}$ norms, so another application of Schauder estimates gives $C^{3,\kappa}$ bounds for $v$. Continuing on nested balls between $B_{2\tau}$ and $B_\tau$ yields \eqref{eq:noncritical}, with constants depending only on $p,k,\delta$; see \cite{GilbargTrudinger}.
\end{proof}

\subsection{High-order growth at critical points}

By \eqref{eq:first} and the Arzelà--Ascoli theorem, the family
\begin{equation}\label{eq:compact}
	\Kcal\coloneqq \overline{\left\{ v|_{\overline B_{1/2}}: \Delta_pv=0\text{ in }B_1,\quad \norm v_{L^\infty(B_1)}\le1,\quad v(0)=0 \right\}}^{\,C^1(\overline B_{1/2})}
\end{equation}
is compact. Every $h\in\Kcal$ satisfies $h(0)=0$ and is $p$-harmonic in $B_{1/2}$: uniform convergence of gradients permits passage to the weak equation because the flux $\xi\mapsto\abs{\xi}^{p-2}\xi$ is continuous for $p>1$.

\begin{lemma}\label{lem:scales}
Let $0<\beta<1$ satisfy $k+\beta<s(p)$. There are $\delta>0$ and finitely many radii $\rho_1,\ldots,\rho_N\in(0,1/4)$, depending only on $p,k,\beta$, such that every $p$-harmonic function $v$ in $B_1$ satisfying
\begin{equation}\label{eq:normalized}
	\norm v_{L^\infty(B_1)}\le1,\qquad v(0)=0,\qquad \abs{Dv(0)}\le\delta
\end{equation}
obeys
\begin{equation}\label{eq:decay}
	\norm v_{L^\infty(B_{\rho_i})}<\rho_i^{k+\beta}
\end{equation}
for at least one $i\in\{1,\ldots,N\}$.
\end{lemma}

\begin{proof}
Set $\sigma'\coloneqq\min\{s(p),k+1\}$, so that $k+\beta<\sigma'\le k+1$. By Iwaniec and Manfredi \cite[Theorem~1]{IwaniecManfredi1989}, every member of the compact set $\Kcal_c\coloneqq\{h\in\Kcal:Dh(0)=0\}$ is of class $C^{k,\sigma'-k}$ near zero. Since $h(0)=0$ and $Dh(0)=0$, \Cref{lem:rigidity} and Taylor's theorem give $h(x)=O(\abs x^{\sigma'})$, with a constant depending on $h$. Thus, for each $h\in\Kcal_c$, there is $\rho_h\in(0,1/4)$ such that
\begin{equation*}
	\norm h_{L^\infty(B_{\rho_h})}<\tfrac12\rho_h^{k+\beta}.
\end{equation*}
The relatively open sets
\begin{equation*}
	U_h\coloneqq\left\{z\in\Kcal:\norm z_{L^\infty(B_{\rho_h})}<\rho_h^{k+\beta}\right\}
\end{equation*}
cover $\Kcal_c$. Choose a finite subcover $U_{h_1},\ldots,U_{h_N}$, set $\rho_i=\rho_{h_i}$, and let $U$ be their union.

Compactness also gives $\delta>0$ such that
\begin{equation*}
	\{z\in\Kcal:\abs{Dz(0)}\le\delta\}\subset U.
\end{equation*}
Otherwise, a sequence in $\Kcal\setminus U$ with gradients at zero tending to zero would converge, along a subsequence, to a member of $\Kcal_c\setminus U$. This proves \eqref{eq:decay}.
\end{proof}

\section{Proof of the main theorem}\label{sec:estimates}

Throughout this section, we fix $k\ge2$ and $0<\beta<1$ with $k+\beta<s(p)$, and write $\sigma=k+\beta$. Let $\delta$ and $\rho_1,\ldots,\rho_N$ be given by \Cref{lem:scales}, and set $\theta=\min_i\rho_i>0$. The qualitative theorem of Iwaniec and Manfredi ensures that $p$-harmonic functions are locally $C^k$.

\begin{proposition}\label{prop:approx}
Let $v$ be a $p$-harmonic function in $B_1$ with $v(0)=0$ and $\norm v_{L^\infty(B_1)}\le1$. Then, the Taylor polynomial $P$ of degree $k$ of $v$ at zero satisfies
\begin{equation}\label{eq:approximation}
	\sum_{j=1}^k\abs{D^jP(0)}\le C,\qquad \norm{v-P}_{L^\infty(B_r)}\le Cr^\sigma\quad \text{for} \quad 0<r\le1,
\end{equation}
where $C$ depends only on $p,k,\beta$.
\end{proposition}

\begin{proof}
We construct the radii inductively, preserving the normalization required by \Cref{lem:scales}. Set $R_0=1$ and $v_0=v$. Suppose that $R_\ell$ has been chosen so that
\begin{equation*}
	v_\ell(x)\coloneqq\frac{v(R_\ell x)}{R_\ell^\sigma}
\end{equation*}
is $p$-harmonic in $B_1$, satisfies $v_\ell(0)=0$, and has $\norm{v_\ell}_{L^\infty(B_1)}\le1$. These properties hold for $\ell=0$ by our assumptions on $v$.

If $\abs{Dv_\ell(0)}>\delta$, we stop. Otherwise, all hypotheses of \Cref{lem:scales} are satisfied, so there is a radius $\rho_i$ such that $\norm{v_\ell}_{L^\infty(B_{\rho_i})}\le\rho_i^\sigma$. Set $R_{\ell+1}=R_\ell\rho_i$. The next rescaling satisfies
\begin{equation*}
	v_{\ell+1}(x)=\frac{v_\ell(\rho_i x)}{\rho_i^\sigma},
\end{equation*}
and is therefore $p$-harmonic in $B_1$ and vanishes at zero. Moreover,
\begin{equation*}
	\norm{v_{\ell+1}}_{L^\infty(B_1)}=\rho_i^{-\sigma}\norm{v_\ell}_{L^\infty(B_{\rho_i})}\le1.
\end{equation*}
Thus the normalization is preserved, and we can repeat the procedure.

At every constructed scale,
\begin{equation}\label{eq:discrete}
	\norm v_{L^\infty(B_{R_\ell})}\le R_\ell^\sigma,\qquad \theta R_\ell\le R_{\ell+1}\le\tfrac14R_\ell,
\end{equation}
where the bounds involving $R_{\ell+1}$ apply whenever the next scale is constructed. Since $Dv_\ell(0)=R_\ell^{1-\sigma}Dv(0)$, the iteration continues indefinitely if $Dv(0)=0$. If $Dv(0)\ne0$, an infinite iteration would give $R_\ell\to0$ and hence $\abs{Dv_\ell(0)}\to\infty$, since $\sigma>1$, contradicting the condition for continuing. Thus the iteration stops after finitely many steps in this case.

If the iteration continues indefinitely, then $R_\ell\to0$. For $R_{\ell+1}\le r\le R_\ell$,
\begin{equation*}
	\norm v_{L^\infty(B_r)}\le R_\ell^\sigma\le\theta^{-\sigma}r^\sigma.
\end{equation*}
Since $v$ is $C^k$ near zero and $\sigma>k$, the Taylor polynomial is $P=0$, so \eqref{eq:approximation} follows.

Otherwise, let $L\ge0$ be the first index with $\abs{Dv_L(0)}>\delta$. By \Cref{lem:noncritical}, the Taylor polynomial $Q$ of degree $k$ of $v_L$ at zero satisfies
\begin{equation*}
	Q(0)=0,\qquad \sum_{j=1}^k\abs{D^jQ(0)}\le C_\delta,\qquad \norm{v_L-Q}_{L^\infty(B_t)}\le C_\delta t^{k+1}\quad \text{for}\quad 0<t\le\tau.
\end{equation*}
Thus $P(x)=R_L^{\sigma}Q(x/R_L)$ is the Taylor polynomial of $v$, and
\begin{equation}\label{eq:coefficients}
	\abs{D^jP(0)}\le C_\delta R_L^{\sigma-j}\le C_\delta\quad \text{for} \quad  j \in \{1,\ldots,k \}.
\end{equation}
For $0<r\le\tau R_L$, the Taylor remainder gives
\begin{equation}\label{eq:small}
	\norm{v-P}_{L^\infty(B_r)}\le C_\delta R_L^{\sigma-k-1}r^{k+1}=C_\delta r^\sigma(r/R_L)^{1-\beta}\le C_\delta r^\sigma.
\end{equation}

For $\tau R_L\le r\le1$, \eqref{eq:discrete} and monotonicity of the supremum give
\begin{equation*}
	\norm v_{L^\infty(B_r)}\le C\max\{r,R_L\}^\sigma.
\end{equation*}
Indeed, use the bound at $R_L$ when $r\le R_L$, and compare consecutive constructed scales when $R_L<r\le1$. Since $\sigma-j>0$, \eqref{eq:coefficients} also yields
\begin{equation*}
	\norm P_{L^\infty(B_r)}\le C\sum_{j=1}^kR_L^{\sigma-j}r^j\le C\max\{r,R_L\}^\sigma.
\end{equation*}
Now $\max\{r,R_L\}\le r/\tau$, which completes the proof of \eqref{eq:approximation}.
\end{proof}

\begin{proof}[Proof of \Cref{thm:main}]
Let $u\in L^p(B_1)$ be $p$-harmonic in $B_1$. Local boundedness, applied to $u$ and $-u$ with integrability exponent $p$ on the concentric balls $B_{3/4}\subset B_{7/8}$ \cite[Lemma~3.6]{Lindqvist}, gives
\begin{equation*}
	\norm u_{L^\infty(B_{3/4})}\le C_b\norm u_{L^p(B_1)},\qquad C_b=C_b(p).
\end{equation*}
Set $M_0=2C_b\norm u_{L^p(B_1)}$; the case $M_0=0$ is immediate. For every $x\in B_{1/2}$, the function
\begin{equation*}
	w_x(y)=\frac{u(x+y/4)-u(x)}{M_0},\qquad y\in B_1,
\end{equation*}
is $p$-harmonic in $B_1$ and vanishes at zero. Since $x+B_{1/4}\subset B_{3/4}$, we also have
\begin{equation*}
	\norm{w_x}_{L^\infty(B_1)}\le\frac{2\norm u_{L^\infty(B_{3/4})}}{M_0}\le1.
\end{equation*}
We can therefore apply \Cref{prop:approx} to the Taylor polynomial $Q_x$ of $w_x$ at zero. Returning to the original variables, set
\begin{equation*}
	P_x(z)=u(x)+M_0Q_x\bigl(4(z-x)\bigr).
\end{equation*}
Then $P_x$ is the Taylor polynomial of $u$ at $x$, and the fixed spatial scaling gives
\begin{align}
	\sum_{j=1}^k\abs{D^ju(x)}&\le CM_0,\label{eq:jet-bound}\\
	\abs{u(z)-P_x(z)}&\le CM_0\abs{z-x}^\sigma\qquad(z\in B_{1/4}(x)).\label{eq:jet-error}
\end{align}

For distinct $x,y\in B_{1/2}$ with $d=\abs{x-y}<1/16$, we have $B_d(x)\subset B_{1/4}(x)\cap B_{1/4}(y)$, since $\abs{z-y}<2d<1/4$ for $z\in B_d(x)$. Both approximation bounds therefore hold on $B_d(x)$, so
\begin{equation*}
	\norm{P_x-P_y}_{L^\infty(B_d(x))}\le CM_0d^\sigma.
\end{equation*}
Since $D^kP_x=D^ku(x)$ and $D^kP_y=D^ku(y)$, scaling and equivalence of norms for polynomials of degree at most $k$ yield
\begin{equation*}
	\abs{D^ku(x)-D^ku(y)}\le Cd^{-k}\norm{P_x-P_y}_{L^\infty(B_d(x))}\le CM_0d^\beta.
\end{equation*}
For $d\ge1/16$, the same H\"older bound follows from \eqref{eq:jet-bound}. Together with local boundedness and the derivative bounds in \eqref{eq:jet-bound}, this proves \eqref{eq:main}. All constants depend only on $p,k,\beta$, as required.
\end{proof}

\section{Consequences}\label{sec:consequences}

In this section, we derive three consequences of \Cref{thm:main}: uniform decay estimates at critical points, a Liouville theorem for entire solutions, and uniform higher-order estimates for the flux of planar $q$-harmonic functions with $q>2$.

\subsection{Decay at critical points}

\begin{corollary}\label{cor:critical}
Let $1<p<2$, $0<\sigma<s(p)$, and let $u\in L^p(B_1)$ be $p$-harmonic in $B_1$. If $x_0\in B_{1/2}$ satisfies $Du(x_0)=0$, then there exists $C=C(p,\sigma)>0$ such that
\begin{equation}\label{eq:critical-decay}
	\abs{u(x)-u(x_0)}\le C\norm{u}_{L^p(B_1)}\abs{x-x_0}^\sigma\qquad\text{for every }x\in B_{1/4}(x_0).
\end{equation}
\end{corollary}

\begin{proof}
Choose an integer $k\ge2$ and $0<\beta<1$ such that $\sigma<k+\beta<s(p)$. By \Cref{thm:main}, $u$ is $C^k$ near $x_0$, and \Cref{lem:rigidity} gives $D^ju(x_0)=0$ for $1\le j\le k$. Thus the Taylor polynomial of degree $k$ at $x_0$ is the constant $u(x_0)$. Applying \eqref{eq:jet-error}, we obtain
\[
	\abs{u(x)-u(x_0)}\le C\norm{u}_{L^p(B_1)}\abs{x-x_0}^{k+\beta}
	\qquad\text{for every }x\in B_{1/4}(x_0).
\]
Since $\abs{x-x_0}<1$ and $k+\beta>\sigma$, this implies \eqref{eq:critical-decay}.
\end{proof}

\subsection{A Liouville theorem}

The conclusion below also follows from the Liouville theorem stated in \cite[Theorem~2.5]{KilpelainenShahgholianZhong2007}, through the classical conjugacy between planar $p$-harmonic and $p'$-harmonic functions, where $p'=p/(p-1)$; see \cite{AronssonLindqvist1988}. We give an independent proof based on \Cref{thm:main}.

\begin{corollary}\label{cor:liouville}
Let $1<p<2$, and let $u$ be $p$-harmonic in $\R^2$. Suppose that there are $\gamma<s(p)$ and $C_1<\infty$ such that
\begin{equation*}
	\sup_{B_R}\abs u\le C_1R^\gamma\qquad\text{for every }R\ge1.
\end{equation*}
Then $u$ is affine.
\end{corollary}

\begin{proof}
Choose an integer $k\ge2$ and $0<\beta<1$ such that $\gamma<k+\beta<s(p)$. For $R\ge1$, the function $u_R(x)=u(Rx)$ is $p$-harmonic in $B_1$, and $\norm{u_R}_{L^p(B_1)}\le\abs{B_1}^{1/p}C_1R^\gamma$. Since $D^ku_R(x)=R^kD^ku(Rx)$, \Cref{thm:main} gives
\begin{equation*}
	[D^ku]_{C^{0,\beta}(B_{R/2})}=R^{-(k+\beta)}[D^ku_R]_{C^{0,\beta}(B_{1/2})}\le CR^{\gamma-k-\beta}.
\end{equation*}
Letting $R\to\infty$, we find that $D^ku$ is constant in $\R^2$, so $u$ is a polynomial of degree $d\le k$.

Suppose that $d\ge2$, and let $P_d\not\equiv0$ be the homogeneous part of degree $d$ of $u$. Then $R^{-d}u(Rx)\to P_d(x)$ in $C^1_{\mathrm{loc}}(\R^2)$ as $R\to\infty$. Each function $x\mapsto R^{-d}u(Rx)$ is $p$-harmonic in $\R^2$, so, as in the proof of \Cref{lem:rigidity}, $P_d$ is $p$-harmonic. This contradicts Lewis's theorem \cite[Lemma~2]{Lewis1980}. Hence $d\le1$.
\end{proof}

\begin{remark}\label{rem:sharp}
The exponent $s(p)$ in \Cref{cor:liouville} is optimal. Indeed, for every $1<p<2$, there exists a non-affine $p$-harmonic function in $\R^2$ of the form
\[
	u(r\cos\vartheta,r\sin\vartheta)=r^{s(p)}f(\vartheta),
\]
where $f$ is a smooth $\pi$-periodic function; see \cite{Aronsson1988} and \cite[Propositions~2.4--2.5]{Tkachev2006}. By homogeneity,
\[
	\sup_{B_R}\abs{u}=R^{s(p)}\sup_{B_1}\abs{u}\qquad\text{for every }R>0,
\]
so the conclusion of \Cref{cor:liouville} fails when $\gamma=s(p)$.
\end{remark}

\subsection{Flux estimates for \texorpdfstring{$q>2$}{q>2}}

The classical stream-function correspondence transfers \Cref{thm:main} to uniform higher-order estimates for the flux of planar $q$-harmonic functions with $q>2$. Write $A_q(\xi)=\abs{\xi}^{q-2}\xi$ and $q'=q/(q-1)\in(1,2)$. Since $1/(q'-1)=q-1$, formula \eqref{eq:IM} gives
\begin{equation}\label{eq:dual-order}
	s(q')-1=\frac16\left(q+\sqrt{q^2+12q-12}\right)>1.
\end{equation}

\begin{corollary}\label{cor:flux}
Let $q>2$, let $j\ge1$ be an integer, and let $0<\beta<1$ satisfy $j+\beta<\frac16\bigl(q+\sqrt{q^2+12q-12}\bigr)$. There exists $C=C(q,j,\beta)>0$ such that every $q$-harmonic function $v\in W^{1,q}(B_1)$ satisfies
\begin{equation}\label{eq:flux}
	\sum_{i=0}^{j}\norm{D^i\bigl(A_q(Dv)\bigr)}_{L^\infty(B_{1/2})}+\bigl[D^j\bigl(A_q(Dv)\bigr)\bigr]_{C^{0,\beta}(B_{1/2})}\le C\norm{Dv}_{L^q(B_1)}^{q-1}.
\end{equation}
\end{corollary}

\begin{proof}
Let $J(\xi_1,\xi_2)=(-\xi_2,\xi_1)$. The field $A_q(Dv)\in L^{q'}(B_1;\R^2)$ is weakly divergence-free and $B_1$ is simply connected, so there exists $u\in W^{1,q'}(B_1)$ such that
\[
	Du=JA_q(Dv)\quad\text{a.e. in }B_1,\qquad\int_{B_1}u\,dx=0.
\]
This is the stream function of $v$; see \cite{AronssonLindqvist1988}. Since $\abs{Du}=\abs{Dv}^{q-1}$ and $(q-1)(q'-2)=2-q$, we have
\[
	\abs{Du}^{q'-2}Du=JDv\quad\text{a.e. in }B_1.
\]
As $\diver(JDv)=0$ in the sense of distributions, $u$ is $q'$-harmonic in $B_1$. Moreover, $\abs{Du}^{q'}=\abs{Dv}^{q}$, so the Poincar\'e inequality gives
\[
	\norm{u}_{L^{q'}(B_1)}\le C\norm{Du}_{L^{q'}(B_1)}=C\norm{Dv}_{L^q(B_1)}^{q-1}.
\]
By \eqref{eq:dual-order}, $j+1+\beta<s(q')$. Applying \Cref{thm:main} to $u$ with exponent $q'$ and $k=j+1$, and using $A_q(Dv)=-JDu$, yields \eqref{eq:flux}.
\end{proof}

\begin{remark}\label{rem:flux}
For $q\ge2$, Balci, Diening, and Weimar \cite[Theorem~2.2]{BalciDieningWeimar2020} proved decay of the oscillation of the flux of planar $q$-harmonic functions with every exponent below one. Qualitative higher-order regularity below the threshold \eqref{eq:dual-order} follows from Iwaniec--Manfredi regularity and classical planar conjugacy. \Cref{cor:flux} provides the corresponding uniform bounds in terms of $\norm{Dv}_{L^q(B_1)}^{q-1}$. This threshold cannot be exceeded, since an improvement would, by conjugacy, improve the optimal regularity order $s(q')$ for planar $q'$-harmonic functions.
\end{remark}

\bigskip

{\small \noindent{\bf Acknowledgments.} JM was partially supported by a gift from the Simons Foundation (Gift Number 962828). This publication is based upon work supported by King Abdullah University of Science and Technology (KAUST) under Award No. ORFS-CRG12-2024-6430. AS acknowledges the hospitality of the mathematics departments at the University of Pittsburgh and Oklahoma State University, where part of this research was conducted. ET acknowledges partial support from the Grayce B. Kerr Chair funds at Oklahoma State University.}

\bigskip

{\small \noindent{\bf AI disclosure.} The authors developed the mathematical ideas, statements, and proofs for \Cref{thm:main} and the results in \Cref{sec:critical,sec:estimates}. AI tools were used to correct typographical errors, improve the writing, review the manuscript, and produce a Lean~4 formalization of \Cref{thm:main} using Mathlib. The consequences in \Cref{sec:consequences} and \Cref{rem:scope} were suggested during an AI-assisted review and verified by the authors. The Lean code and a description of the formalization process are available at \url{https://github.com/na0b/phigherorder-lean}.}

\end{document}